\documentclass[english]{article}
\usepackage[latin1]{inputenc}
\usepackage{amsmath,amsthm,amssymb,babel}
\newtheorem{teo}{Theorem}

\newtheorem{lemma}{Lemma}

\def\proof{{\it Proof.}\ }
\def\endproof{\hfill $\Box$\par\vskip3mm}

\def\neweq#1{\begin{equation}\label{#1}}
\def\endeq{\end{equation}}

\def\phi{\varphi}
\def\RR{{\mathbb R} }

\date{}

\title{\sc On a PDE involving the ${\cal A}_{p(\cdot)}$-Laplace operator}
\author{\sc Mihai Mih\u ailescu$\,^{a}$ and Du\v san Repov\v s$\,^{b}$\\
\small
$^a\,$Department of Mathematics, University of Craiova,  200585 Craiova,
Romania\\
\small
$^b\,$Faculty of Mathematics and Physics, and Faculty of Education, University of Ljubljana,\\
\small
 POB 2964, Ljubljana, Slovenia 1001\\
\small
E-mail addresses:  {\tt mmihailes@yahoo.com}\qquad {\tt dusan.repovs@guest.arnes.si} }

\begin{document}
\maketitle \noindent{\small{\sc Abstract}.  This paper establishes existence of solutions for a partial differential equation in which a differential operator involving variable exponent growth conditions is present. This operator represents a generalization of the $p(\cdot)$-Laplace operator, i.e. $\Delta_{p(\cdot)}u={\rm div}(|\nabla u|^{p(\cdot)-2}\nabla u)$, where $p(\cdot)$ is a continuous function. The proof of the main result is based on Schauder's fixed point theorem combined with adequate variational arguments. The function space setting used here makes appeal to the variable exponent Lebesgue and Sobolev spaces.
}

\section{Introduction}\label{sec1}
Let $\Omega\subset{\mathbb R}^N$ ($N\geq 2$) be a bounded domain with smooth boundary $\partial\Omega$. Let $p(\cdot):\overline\Omega\rightarrow(2,\infty)$ be a continuous function such that
\begin{equation}\label{lambda1}
\lambda_1:=\inf_{u\in C_0^\infty(\Omega)\setminus\{0\}}\frac{\displaystyle\int_\Omega|\nabla u|^{p(x)}\;dx}{\displaystyle\int_\Omega|u|^{p(x)}\;dx}>0\,.
\end{equation}
We point out that  property $\lambda_1>0$ is not true for all functions $p(\cdot)$. For instance, assuming that there exists an open set $U\subset\Omega$ and a point $x_0\in U$ such that $p(x_0)<p(x)$ (or $p(x_0)>p(x)$) for all $x\in\partial U$, then by \cite[Theorem 3.1]{FZZ}  we get $\lambda_1=0$. On the other hand, there are results establishing sufficient conditions on $p(\cdot)$ in order to satisfy $\lambda_1>0$. Indeed,  it was proved in \cite[Theorem 3.3]{FZZ} that assuming that  there exists a vector $l\in\RR^N\setminus\{0\}$ such that, for any $x\in\Omega$, the function $f(t)=p(x+tl)$ is monotone, for $t\in I_x:=\{s;\;x+sl\in\Omega\}$ then $\lambda_1>0$. Furthermore,  it was shown in \cite[Theorem 1]{CVEE} that (\ref{lambda1}) holds provided that $p(\cdot)\in C^1(\Omega;{\mathbb R})$ and that there exists $\overrightarrow a\in C^1(\Omega;{\mathbb R}^N)$ such that ${\rm div}\overrightarrow a(x)\geq a_0>0$ and $\overrightarrow a(x)\cdot\nabla p(x)=0$, for every $x\in\Omega$ (see also \cite[Theorem 1]{nonlinearity} for similar results). Finally, we recall a very well-known fact that in the special case when $p(\cdot)$ is a constant function (defined on the interval $(1,\infty)$) then (\ref{lambda1}) holds.

Next, assume that $A:\Omega\rightarrow{\mathbb R}^{N^2}$ is a symmetric function matrix, i.e. $a_{ij}=a_{ji}$, such that $a_{ij}\in L^\infty(\Omega)\cap C^1(\Omega)$ and
\begin{equation}\label{eliptic}
\langle A\xi,\xi\rangle=\sum_{i,j=1}^Na_{ij}(x)\xi_i\xi_j\geq|\xi|^2,\;\;\;\forall\;x\in\Omega,\;\xi\in{\mathbb R}^N\,,
\end{equation}
where $\langle\cdot,\cdot\rangle$ denotes the scalar product on ${\mathbb R}^N$.

In this paper we are concerned with the study of nonlinear and nonhomogeneous problems of type
\begin{equation}\label{problem}
\left\{\begin{array}{ll} -{\rm div}(\alpha(x,u)\langle A\nabla u,\nabla u\rangle^{\frac{p(x)-2}{2}}A\nabla u)=f &\mbox{for
}\phantom{\partial}x\in\Omega,\\
u=0 &\mbox{for}\phantom{\partial}x\in\partial\Omega\,,
\end{array}\right.
\end{equation}
where $\alpha:\Omega\times\mathbb R\rightarrow(0,\infty)$ is a bounded function and $f:\Omega\rightarrow\mathbb R$ is a measurable function belonging to a suitable  Lebesgue type space which  will be specified later on in the paper. The differential operator involved in equation (\ref{problem}) will be denoted by ${\cal A}_{p(\cdot)}:={\rm div}(\alpha(x,u)\langle A\nabla u,\nabla u\rangle^{\frac{p(x)-2}{2}}A\nabla u)$ and will be called the {\it ${\cal A}_{p(\cdot)}$-Laplace operator}. It  represents a generalization of the {\it $p(\cdot)$-Laplace operator}, i.e. $\Delta_{p(\cdot)}u={\rm div}(|\nabla u|^{p(\cdot)-2}\nabla u)$, which is obtained in the case when $A=Id$ and $\alpha\equiv 1$.  In the last decades special attention has been paid to $p(\cdot)$-Laplace type operators  since they can model with sufficient accuracy the phenomena arising from the study of electrorheological fluids (Ru\v{z}i\v{c}ka \cite{R}, Rajagopal \& Ru\v{z}i\v{c}ka \cite{RaRu01}), image restoration (Chen {\it et al.} \cite{CLR}), mathematical biology (Fragnelli \cite{frag}), dielectric breakdown, electrical resistivity and polycrystal plasticity (Bocea \& Mih\u ailescu \cite{BM-NA2010}, Bocea {\it et al.} \cite{B-M-P: Ricerche2010})  or they  arise in the study of some models for growth of heterogeneous sandpiles (Bocea {\it et al.} \cite{BMPLR}). In a similar   context,   we note that a collection of  results obtained in the field of partial differential equations involving $p(\cdot)$-Laplace type operators can be found in the survey paper by Harjulehto {\it et al.} \cite{HHLN}. Finally, we recall that in the case  when $p(\cdot)$ is a constant function, problems involving ${\cal A}_{p}$-Laplace type operators have been widely studied. In this regard  we point out the papers by Reshetnyak \cite{res}, Alvino {\it et al.} \cite{AFT} and El Khalil {\it et al.} \cite{elk} and the references therein.

\section{A review on variable exponent spaces}\label{sec2}
In this section we provide a brief review of  basic properties of
the variable exponent Lebesgue-Sobolev spaces. For more details we
refer to the book by  Diening {\it et al.} \cite{book} and the paper by  Kovacik and  R\'akosn\'{\i}k
\cite{KR}.
\smallskip

In this paper we reduce all our discussion to the special  case when $\Omega\subset\RR^N$ is  an open bounded set.  For any  continuous function
$p:\overline\Omega\rightarrow(1,\infty)$ we define
$$p^-:=\inf\limits_{x\in\Omega}p(x)\;\;\;{\rm and}\;\;\;
p^+:=\sup\limits_{x\in\Omega}p(x).$$ Next, we define the variable exponent Lebesgue space
$L^{p(\cdot)}(\Omega)$ by
$$L^{p(\cdot)}(\Omega)=\left\{u:\Omega \to \mathbb{R}\ {\rm measurable}  \ : \
\int_\Omega|u(x)|^{p(x)}\;dx<\infty\right\}\,.$$
Clearly, $L^{p(\cdot)}(\Omega)$ is a Banach space when endowed with the
so-called {\it Luxemburg norm}, defined by
$$|u|_{p(\cdot)}:=\inf\left\{\mu>0 \ : \ \int_\Omega\left|
\frac{u(x)}{\mu}\right|^{p(x)}\;dx\leq 1\right\}.$$ We note that the
variable exponent Lebesgue space is a special case of an
Orlicz-Musielak space. For constant functions $p,\
L^{p(\cdot)}(\Omega)$ reduces to the classical Lebesgue space
$L^{p}(\Omega),$ endowed with the standard norm
$$\|u\|_{L^{p}(\Omega)}:=\left(\int_\Omega |u(x)|^{p}dx\right)^{1/p}.$$

We recall that $L^{p(\cdot)}(\Omega)$ is separable and reflexive. Since $\Omega$ is bounded, if $p_1$, $p_2$ are variable exponents such 
that $p_1 \leq p_2$ in $\Omega,$  the embedding
$L^{p_2(\cdot)}(\Omega)\hookrightarrow L^{p_1(\cdot)}(\Omega)$ is
continuous and its norm does not exceed $|\Omega|+1$.

We denote by $L^{p^{'}(\cdot)}(\Omega)$ the conjugate space of
$L^{p(\cdot)}(\Omega)$, where $1/p(x)+1/p^{'}(x)=1$. For any $u\in
L^{p(\cdot)}(\Omega)$ and $v\in L^{p^{'}(\cdot)}(\Omega)$ the following  H\"older
type inequality
\begin{equation}\label{Hol}
\left|\int_\Omega u v\;dx\right|\leq\left(\frac{1}{p^-}+
\frac{1}{{p^{'}}^-}\right)|u|_{p(\cdot)}|v|_{p^{'}(\cdot)}
\end{equation}
holds.

A key role in manipulating the variable exponent Lebesgue and
Sobolev (see below) spaces is played by the {\it modular} of the
space $L^{p(\cdot)}(\Omega)$, which is the mapping
 $\rho_{p(\cdot)}:L^{p(\cdot)}(\Omega)\rightarrow\RR$ defined by
$$\rho_{p(\cdot)}(u):=\int_\Omega|u(x)|^{p(x)}\;dx.$$

If $u\in L^{p(\cdot)}(\Omega)$ then the
following relations hold:
\begin{equation}\label{L4}
|u|_{p(\cdot)}>1\;\;\;\Rightarrow\;\;\;|u|_{p(\cdot)}^{p^-}\leq\rho_{p(\cdot)}(u)
\leq|u|_{p(\cdot)}^{p^+}\,;
\end{equation}
\begin{equation}\label{L5}
|u|_{p(\cdot)}<1\;\;\;\Rightarrow\;\;\;|u|_{p(\cdot)}^{p^+}\leq
\rho_{p(\cdot)}(u)\leq|u|_{p(\cdot)}^{p^-}\,;
\end{equation}
\begin{equation}\label{L6}
|u|_{p(\cdot)}=1\Leftrightarrow\rho_{p(\cdot)}(u)=1\,.
\end{equation}
The variable exponent Sobolev space $W^{1,p(\cdot)}(\Omega)$ is
defined by
$$
W^{1,p(\cdot)}(\Omega):=\{u\in L^{p(\cdot)}(\Omega):\;|\nabla u|\in
L^{p(\cdot)} (\Omega) \}.
$$
On this space one can consider  the following  norm
$$
\|u\|_{p(\cdot)}: =|u|_{p(\cdot)}+|\nabla u|_{p(\cdot)},
$$
where, in the above definition  $\ |\nabla
u|_{p(\cdot)}$ stands for the Luxemburg norm of  $|\nabla u|$. We
note that in the context of this discussion that 
$W^{1,p(\cdot)}(\Omega)$ is also a separable and reflexive Banach
space.

Finally, we define $W_0^{1,p(\cdot)}(\Omega)$ as the closure of
$C_0^\infty(\Omega)$ under the norm
$$\|u\|=|\nabla u|_{p(\cdot)}.$$
Note that  $(W_0^{1,p(\cdot)}(\Omega),\|\cdot\|)$ is also a separable and
reflexive Banach space. We remark  that if $q:\overline\Omega\rightarrow(1,\infty)$ is a continuous function such that
 $q(x)<p^\star(x)$ for all $x\in\overline\Omega$ then the
embedding
$W_0^{1,p(\cdot)}(\Omega)\hookrightarrow L^{q(\cdot)}(\Omega)$
is compact and continuous, where $p^\star(x)=\frac{Np(x)}{N-p(x)}$
if $p(x)<N$ or $p^\star(x)=+\infty$ if $p(x)\geq N$.

\section{The main result}\label{sec3}
The main result of this paper is given by the following theorem.
\begin{teo}\label{t1}
Assume that  $\alpha:\Omega\times\mathbb R\rightarrow\mathbb R$   is a Carath\'eodory function for which there exist two positive constants $0<\lambda\leq\Lambda$ such that:
\begin{equation}\label{ineq}
0<\lambda\leq\alpha(x,t)\leq\Lambda,\;\;\;{\rm a.e.}\;x\in\Omega,\;\forall\;t\in\mathbb R\,.
\end{equation}
 Assume that conditions (\ref{lambda1}) and (\ref{eliptic}) from Section \ref{sec1} are satisfied. Then for each $f\in L^{p^{'}(\cdot)}(\Omega)$ there exists  a weak solution of problem (\ref{problem}), i.e. a function $u\in W_0^{1,p(\cdot)}(\Omega)$ such that 
$$\int_\Omega\alpha(x,u)\langle A\nabla u,\nabla u\rangle^{\frac{p(x)-2}{2}}\langle A\nabla u,\nabla\phi\rangle\;dx=
\int_\Omega f\phi\;dx\,,$$
for all $\phi\in W_0^{1,p(\cdot)}(\Omega)$.
\end{teo}

\section{Proof of the main result}\label{sec4}
 Fix an arbitrary function $f\in L^{p^{'}(\cdot)}(\Omega)$. The main ingredient  of our proof of  Theorem \ref{t1} will be  Schauder's fixed point theorem:
\bigskip

\noindent{\bf Schauder's Fixed Point Theorem.} {\it Assume that $K$ is a compact and convex subset of the Banach space $B$ and  $S:K\rightarrow K$ is a continuous map. Then $S$ possesses a fixed point.}
\bigskip

We start by proving some auxiliary results which will  be useful in establishing Theorem \ref{t1}.

\begin{lemma}\label{l1}
For each $v\in L^{p(\cdot)}(\Omega)$ the problem
\begin{equation}\label{problem2}
\left\{\begin{array}{ll} -{\rm div}(\alpha(x,v)\langle A\nabla u,\nabla u\rangle^{\frac{p(x)-2}{2}}A\nabla u)=f &\mbox{for
}\phantom{\partial}x\in\Omega,\\
u=0 &\mbox{for}\phantom{\partial}x\in\partial\Omega\,,
\end{array}\right.
\end{equation}
 has a weak solution $u\in W_0^{1,p(\cdot)}(\Omega)$, i.e.
 \begin{equation}\label{4}
 \int_\Omega\alpha(x,v)\langle A\nabla u,\nabla u\rangle^{\frac{p(x)-2}{2}}\langle A\nabla u,\nabla\phi\rangle\;dx=
\int_\Omega f\phi\;dx\,,
 \end{equation}
 for all $\phi\in W_0^{1,p(\cdot)}(\Omega)$.
\end{lemma}
\proof
Fix $v\in L^{p(\cdot)}(\Omega)$. First, we note that condition (\ref{ineq}) from Theorem \ref{t1} guarantees  that $\alpha(x,v)\in L^\infty(\Omega)$.

Consider the energy functional associated with problem (\ref{problem2}), $J:W_0^{1,p(\cdot)}(\Omega)\rightarrow\mathbb R$,
$$J(u)=\int_\Omega\frac{\alpha(x,v)}{p(x)}\langle A\nabla u,\nabla u\rangle^{p(x)/2}\;dx-\int_\Omega fu\;dx\,.$$
Standard arguments imply that $J\in C^1(W_0^{1,p(\cdot)}(\Omega);\mathbb R)$ with the derivative given by
$$\langle J^{'}(u),\phi\rangle= \int_\Omega\alpha(x,v)\langle A\nabla u,\nabla u\rangle^{\frac{p(x)-2}{2}}\langle A\nabla u,\nabla\phi\rangle\;dx=
\int_\Omega f\phi\;dx\,,
$$
for all $u,\; \phi\in W_0^{1,p(\cdot)}(\Omega)$. Thus, weak solutions of problem (\ref{problem2}) are exactly the critical points of the functional $J$.

Since (\ref{eliptic}) and (\ref{ineq}) are fulfilled it follows that for each $u\in W_0^{1,p(\cdot)}(\Omega)$ with $\|u\|>1$ we have
\begin{eqnarray*}
J(u)&\geq&\frac{\lambda}{p^+}\int_\Omega|\nabla u|^{p(x)}\;dx-\int_\Omega fu\;dx\\
&\geq&\frac{\lambda}{p^+}\|u\|^{p^-}-c|f|_{p^{'}(\cdot)}\|u\|\,,
\end{eqnarray*}
where $c$ is a positive constant. The above estimate shows that $J$ is coercive.

On the other hand,  it was pointed out in \cite[p. 449]{AFT} that the following Clarkson's type inequality 
\begin{equation}\label{convexity}
\frac{\langle A\xi_1,\xi_1\rangle^{s/2}+\langle A\xi_2,\xi_2\rangle^{s/2}}{2}\geq\left\langle A\left(\frac{\xi_1+\xi_2}{2}\right),\frac{\xi_1+\xi_2}{2}\right\rangle^{s/2}+\left\langle A\left(\frac{\xi_1-\xi_2}{2}\right),\frac{\xi_1-\xi_2}{2}\right\rangle^{s/2}\,,
\end{equation}
holds for all $s\geq 2$ and $\xi_1,\;\xi_2\in{\mathbb R}^N$. Thus, we deduce that $J$ is convex and consequently weakly lower semi-continuous.

Since $J$ is coercive and weakly lower semi-continuous we conclude  via the Direct Method of the Calculus of Variations (see, e.g. \cite[Theorem 1.2~]{S}), that there exists  a global minimum point of $J$, $u\in W_0^{1,p(\cdot)}(\Omega)$ and consequently a weak solution of problem (\ref{problem2}). The proof of Lemma \ref{l1} is thus complete.   \endproof

Next, for each $v\in L^{p(\cdot)}(\Omega)$ let $u=T(v)\in W_0^{1,p(\cdot)}(\Omega)$ be the weak solution of problem (\ref{problem2}) given by Lemma \ref{l1}. Thus,  we can actually introduce an application $T:L^{p(\cdot)}(\Omega)\rightarrow W_0^{1,p(\cdot)}(\Omega)$ associating to each  $v\in L^{p(\cdot)}(\Omega)$, the solution of problem (\ref{problem2}), $T(v)\in W_0^{1,p(\cdot)}(\Omega)$.

\begin{lemma}\label{l2}
There exists $C>0$ a universal constant such that
\begin{equation}\label{2stele}
\int_\Omega|\nabla T(v)|^{p(x)}\;dx\leq C,\;\;\;\forall\;v\in L^{p(\cdot)}(\Omega)\,.
\end{equation}
\end{lemma}
\proof
Taking $\phi=T(v)$ in (\ref{4}) we find
$$\int_\Omega\alpha(x,v)\langle A\nabla T(v),\nabla T(v)\rangle^{\frac{p(x)}{2}}\;dx=
\int_\Omega fT(v)\;dx,\;\;\;\forall\;v\in L^{p(\cdot)}(\Omega)\,.$$
Taking into account  relation (\ref{ineq}) and condition (\ref{eliptic}), the above equality yields
\begin{equation}\label{estim}\lambda\int_\Omega|\nabla T(v)|^{p(x)}\;dx\leq\int_\Omega fT(v)\;dx,\;\;\;\forall\;v\in L^{p(\cdot)}(\Omega)\,.\end{equation}
Let now $\epsilon>0$ be such that $\epsilon<\min\{1,\lambda,\lambda/\lambda_1\}$. Then, by Young's inequality (see, e.g. \cite[the footnote on p. 56]{B})  we deduce
$$f(x)T(v(x))\leq\frac{1}{\epsilon^{p(x)-1}}|f(x)|^{p^{'}(x)}+\epsilon|T(v(x))|^{p(x)},\;\;\;\forall\;v\in L^{p(\cdot)}(\Omega),\;x\in\Omega\,,$$
or, since $\epsilon\in(0,1)$ there exists $C_\epsilon:=\frac{1}{\epsilon^{p^+-1}}$ such that
$$f(x)T(v(x))\leq C_\epsilon|f(x)|^{p^{'}(x)}+\epsilon|T(v(x))|^{p(x)},\;\;\;\forall\;v\in L^{p(\cdot)}(\Omega),\;x\in\Omega\,.$$
Integrating  the above estimate over $\Omega$ and taking into account that relations (\ref{lambda1}) and (\ref{estim}) hold  we get
$$\lambda\int_\Omega|\nabla T(v)|^{p(x)}\;dx\leq C_\epsilon\int_\Omega|f|^{p^{'}(x)}\;dx+\frac{\epsilon}{\lambda_1}\int_\Omega|\nabla T(v)|^{p(x)}\;dx,\;\;\;\forall\;v\in L^{p(\cdot)}(\Omega)\,.$$
Consequently, taking
$$C:=\displaystyle\frac{C_\epsilon\displaystyle\int_\Omega|f|^{p^{'}(x)}\;dx}{\lambda-\displaystyle\frac{\epsilon}{\lambda_1}}\,,$$
we infer that relation (\ref{2stele}) holds true. The proof of Lemma \ref{l2} is thus also complete.  \endproof

\noindent{\bf Remark 1.} By Lemma \ref{l2} and relation (\ref{lambda1}) it  clearly follows that there exists a universal constant $C_1>0$ such that
$$\int_\Omega|T(v)|^{p(x)}\;dx\leq C_1,\;\;\;\forall\;v\in L^{p(\cdot)}(\Omega)\,.$$

\begin{lemma}\label{l3}
The map $T:L^{p(\cdot)}(\Omega)\rightarrow W_0^{1,p(\cdot)}(\Omega)$ is continuous.
\end{lemma}
\proof
Let $(v_n),\; v\subset L^{p(\cdot)}(\Omega)$ be  such that $v_n$ converges to $v$ in $L^{p(\cdot)}(\Omega)$ as $n\rightarrow\infty$. Set
$$u_n:= T(v_n),\;\;\;\forall\; n\,.$$
By Lemma \ref{l2} we have
$$\int_\Omega|\nabla u_n|^{p(x)}\;dx=\int_\Omega|\nabla T(v_n)|^{p(x)}\;dx\leq C,\;\;\;\forall\;n\,,$$
i.e. $(u_n)$ is bounded on $W_0^{1,p(\cdot)}(\Omega)$. It follows that by  eventually passing to a subsequence we can conclude  that $u_n$ converges weakly to $u$ in $W_0^{1,p(\cdot)}(\Omega)$.

On the other hand, for each $n$ we have
\begin{equation}\label{163276}
 \int_\Omega\alpha(x,v_n)\langle A\nabla u_n,\nabla u_n\rangle^{\frac{p(x)-2}{2}}\langle A\nabla u_n,\nabla\phi\rangle\;dx=
\int_\Omega f\phi\;dx\,,
 \end{equation}
 for all $\phi\in W_0^{1,p(\cdot)}(\Omega)$. Taking $\phi=u_n-u$ in the above equality it follows that
 $$\int_\Omega\alpha(x,v_n)\langle A\nabla u_n,\nabla u_n\rangle^{\frac{p(x)-2}{2}}\langle A\nabla u_n,\nabla u_n-\nabla u\rangle\;dx=o(1)\,.$$
This fact and relation (\ref{ineq}) yield
\begin{equation}\label{o}
\int_\Omega\langle A\nabla u_n,\nabla u_n\rangle^{\frac{p(x)-2}{2}}\langle A\nabla u_n,\nabla u_n-\nabla u\rangle\;dx=o(1)\,.
\end{equation}
Next, by  taking $\phi=u_n$ in (\ref{163276}) we get
$$ \int_\Omega\alpha(x,v_n)\langle A\nabla u_n,\nabla u_n\rangle^{\frac{p(x)}{2}}\;dx=
\int_\Omega fu_n\;dx\,,$$
for each $n$. Relation (\ref{ineq}), H\"older's inequality, Poincar\'e's inequality and the fact that $(u_n)$ is bounded on $W_0^{1,p(\cdot)}(\Omega)$ imply that there exist some constants $C_2,\; C_3\; C_4>0$ such that
$$\lambda \int_\Omega\langle A\nabla u_n,\nabla u_n\rangle^{\frac{p(x)}{2}}\;dx\leq\int_\Omega fu_n\;dx\leq
C_2|f|_{p^{'}(\cdot)}|u_n|_{p(\cdot)}\leq C_3\|u_n\|\leq C_4,\;\;\;\forall\;n\,.$$
The above estimates assure that  sequence $(\int_\Omega\langle A\nabla u_n,\nabla u_n\rangle^{\frac{p(x)}{2}}\;dx)$ is bounded. Therefore  we can  deduce that there exists $b>0$ such that, up to a subsequence,
$$\lim_{n\rightarrow\infty}\int_\Omega\langle A\nabla u_n,\nabla u_n\rangle^{\frac{p(x)}{2}}\;dx=b\,.$$
Furthermore, recalling relation (\ref{convexity}) and the fact that $p(x)\geq 2$ for all $x\in\Omega$, we deduce that the map
\begin{equation}\label{lsc}W_0^{1,p(\cdot)}(\Omega)\ni w\rightarrow\int_\Omega\langle A\nabla w,\nabla w\rangle^{\frac{p(x)}{2}}\;dx\in{\mathbb R}\,,\end{equation}
is convex and consequently weakly lower semi-continuous. Thus, we deduce
$$\int_\Omega\langle A\nabla u,\nabla u\rangle^{\frac{p(x)}{2}}\;dx\leq\liminf_{n\rightarrow\infty}\int_\Omega\langle A\nabla u_n,\nabla u_n\rangle^{\frac{p(x)}{2}}\;dx=b\,.$$
On the other hand, using  relation (A.2) in \cite{AFT}, i.e.
$$\langle A\xi_2,\xi_2\rangle^s\geq\langle A\xi_1,\xi_1\rangle^s+s\langle A\xi_1,\xi_1\rangle^{s-2}\langle A\xi_1,\xi_2-\xi_1\rangle,\;\;\;\forall\;\xi_1,\xi_2\in{\mathbb R}^n,\;s\geq 2\,,$$
we obtain that
$$\int_\Omega\langle A\nabla u,\nabla u\rangle^{\frac{p(x)}{2}}\;dx\geq \int_\Omega\langle A\nabla u_n,\nabla u_n\rangle^{\frac{p(x)}{2}}\;dx+
p^-\int_\Omega\langle A\nabla u_n,\nabla u_n\rangle^{\frac{p(x)-2}{2}}\langle A\nabla u_n,\nabla u-\nabla u_n\rangle\;dx,\;\;\;\forall\;n\,.$$
The above pieces of information and relation (\ref{o}) show that
$$\int_\Omega\langle A\nabla u,\nabla u\rangle^{\frac{p(x)}{2}}\;dx=b\,.$$
Taking into account that $(\frac{u_n+u}{2})$ converges weakly to $u$ in $W_0^{1,p(\cdot)}(\Omega)$ and  again invoking  the weak lower semi-continuity of the map defined in relation (\ref{lsc}) we find
\begin{equation}\label{9stele}
b=\int_\Omega\langle A\nabla u,\nabla u\rangle^{\frac{p(x)}{2}}\;dx\leq\liminf_{n\rightarrow\infty}\int_\Omega\left\langle A\nabla\frac{u_n+u}{2},\nabla \frac{u_n+u}{2}\right\rangle^{\frac{p(x)}{2}}\;dx\,.
\end{equation}
Assume by contradiction that $(u_n)$ does not converge (strongly) to $u$ in $W_0^{1,p(\cdot)}(\Omega)$. Then there exist $\epsilon>0$ and a subsequence of $(u_n)$, still denoted by $(u_n)$, such that
$$\int_\Omega\left|\frac{\nabla u_n-\nabla u}{2}\right|^{p(x)}\;dx\geq\epsilon,\;\;\;\forall\;n\,.$$
On the other hand, relations (\ref{convexity}) and (\ref{eliptic}) imply
\begin{eqnarray*}
\frac{1}{2}\int_\Omega\langle A\nabla u,\nabla u\rangle^{\frac{p(x)}{2}}\;dx&+&\frac{1}{2}\int_\Omega\langle A\nabla u_n,\nabla u_n\rangle^{\frac{p(x)}{2}}\;dx-\int_\Omega\left\langle A\nabla\frac{u_n+u}{2},\nabla \frac{u_n+u}{2}\right\rangle^{\frac{p(x)}{2}}\;dx\\
&\geq&
\int_\Omega\left\langle A\nabla\frac{u-u_n}{2},\nabla \frac{u-u_n}{2}\right\rangle^{\frac{p(x)}{2}}\;dx\\
&\geq&\int_\Omega\left|\frac{\nabla u-\nabla u_n}{2}\right|^{p(x)}\;dx,\;\;\;\forall\; n\,.
\end{eqnarray*}
The last two estimates yield
$$b-\epsilon\geq\limsup_{n\rightarrow\infty}\int_\Omega\left\langle A\nabla\frac{u_n+u}{2},\nabla \frac{u_n+u}{2}\right\rangle^{\frac{p(x)}{2}}\;dx\,,$$
which contradicts (\ref{9stele}). Consequently, $(u_n)$  converges (strongly) to $u$ in $W_0^{1,p(\cdot)}(\Omega)$, or
$T:L^{p(\cdot)}(\Omega)\rightarrow W_0^{1,p(\cdot)}(\Omega)$ is continuous. The proof of Lemma \ref{l3} is complete.  \endproof

\noindent{\bf Remark 2.} Since $W_0^{1,p(\cdot)}(\Omega)$ is compactly embedded in $L^{p(\cdot)}(\Omega)$ (i.e. the inclusion operator $i:W_0^{1,p(\cdot)}(\Omega)\rightarrow L^{p(\cdot)}(\Omega)$ is compact),  it follows by Lemma \ref{l3} that the operator $S:L^{p(\cdot)}(\Omega)\rightarrow L^{p(\cdot)}(\Omega)$, $S=i\circ T$ is compact.
\bigskip

{\bf Proof of Theorem \ref{t1}.}
Let $C_1$ be the constant given in Remark 1, i.e.
$$\int_\Omega|S(v)|^{p(x)}\;dx\leq C_1,\;\;\;\forall\;v\in L^{p(\cdot)}(\Omega)\,.$$
Consider the ball
$$B_{C_1}(0):=\{v\in L^{p(\cdot)}(\Omega):\;\;\int_\Omega|v|^{p(x)}\;dx\leq C_1\}\,.$$
Clearly, $B_{C_1}(0)$ is a convex  closed subset of $L^{p(\cdot)}(\Omega)$ and $S(B_{C_1}(0))\subset B_{C_1}(0)$. Moreover, by Remark 2,  $S(B_{C_1}(0))$ is relatively compact in $B_{C_1}(0)$.

Finally, by Lemma \ref{l3} and Remark 2, $S:B_{C_1}(0)\rightarrow B_{C_1}(0)$ is a continuous map. Hence we can apply Schauder's fixed point theorem to obtain $S$ with a fixed point. This gives us a weak solution to problem (\ref{problem}) and thus the proof of Theorem \ref{t1} is finally  complete. \endproof

\medskip
{\bf Acknowledgments.} The research was supported by the Slovenian Research grants P1-0292-0101 and J1-4144-0101.

\medskip

\end{document}